\documentclass[11pt]{article}

\usepackage{amssymb}

\newcommand{\newc}{\newcommand}

\newc{\eqnoset}{\setcounter{equation}{0}}

\newcommand{\mref}[1]{(\ref{#1})}
\newcommand{\reflemm}[1]{Lemma~\ref{#1}}
\newcommand{\refrem}[1]{Remark~\ref{#1}}
\newcommand{\reftheo}[1]{Theorem~\ref{#1}}

\newcommand{\refcoro}[1]{Corollary~\ref{#1}}
\newcommand{\refprop}[1]{Proposition~\ref{#1}}
\newcommand{\refsec}[1]{Section~\ref{#1}}

\newcommand{\beq}{\begin{equation}}
\newcommand{\eeq}{\end{equation}}
\newcommand{\beqno}[1]{\begin{equation}\label{#1}}

\newcommand{\barr}{\begin{array}}
\newcommand{\earr}{\end{array}}

\newc{\bearr}{\begin{eqnarray*}}
\newc{\eearr}{\end{eqnarray*}}

\newc{\bearrno}[1]{\begin{eqnarray}\label{#1}}
\newc{\eearrno}{\end{eqnarray}}

\newc{\non}{\nonumber}
\newc{\nol}{\nonumber\nl}

\newcommand{\bdes}{\begin{description}}
\newcommand{\edes}{\end{description}}
\newc{\benu}{\begin{enumerate}}
\newc{\eenu}{\end{enumerate}}
\newc{\btab}{\begin{tabular}}
\newc{\etab}{\end{tabular}}

\newtheorem{theorem}{Theorem}[section]
\newtheorem{defi}[theorem]{Definition}
\newtheorem{lemma}[theorem]{Lemma}
\newtheorem{rem}[theorem]{Remark}
\newtheorem{exam}[theorem]{Example}
\newtheorem{propo}[theorem]{Proposition}
\newtheorem{corol}[theorem]{Corollary}

\newcommand{\btheo}[1]{\begin{theorem}\label{#1}}
\newc{\brem}[1]{\begin{rem}\label{#1}\em}
\newc{\bexam}[1]{\begin{exam}\label{#1}\em}
\newc{\bdefi}[1]{\begin{defi}\label{#1}}
\newcommand{\blemm}[1]{\begin{lemma}\label{#1}}
\newcommand{\bprop}[1]{\begin{propo}\label{#1}}
\newcommand{\bcoro}[1]{\begin{corol}\label{#1}}
\newcommand{\etheo}{\end{theorem}}
\newcommand{\elemm}{\end{lemma}}
\newcommand{\eprop}{\end{propo}}
\newcommand{\ecoro}{\end{corol}}
\newc{\erem}{\end{rem}}
\newc{\eexam}{\end{exam}}
\newc{\edefi}{\end{defi}}

\newc{\rmk}[1]{{\bf REMARK #1: }}
\newc{\DN}[1]{{\bf DEFINITION #1: }}

\newcommand{\bproof}{{\bf Proof:~~}}
\newc{\eproof}{{\vrule height8pt width5pt depth0pt}\vspace{3mm}}

\newc{\bfrac}[2]{\dspl{\frac{#1}{#2}}}

\newc{\nid}{\noindent}

\newcommand{\dspl}{\displaystyle}
\newc{\grad}{\nabla}
\newc{\Div}{\mbox{div}}
\newc{\pdt}[1]{\dspl{\frac{\partial{#1}}{\partial t}}}
\newc{\pdn}[1]{\dspl{\frac{\partial{#1}}{\partial \nu}}}
\newc{\pdNi}[1]{\dspl{\frac{\partial{#1}}{\partial \mathcal{N}_i}}}
\newc{\pD}[2]{\dspl{\frac{\partial{#1}}{\partial #2}}}
\newc{\dt}{\dspl{\frac{d}{dt}}}
\newc{\bdry}[1]{\mbox{$\partial #1$}}
\newc{\sgn}{\mbox{sign}}

\newc{\Hess}[1]{\frac{\partial^2 #1}{\pdh z_i \pdh z_j}}
\newc{\hess}[1]{\partial^2 #1/\pdh z_i \pdh z_j}

\newc{\ag}{\alpha}
\newc{\bg}{\beta}
\newc{\cg}{\gamma}\newc{\Cg}{\Gamma}
\newc{\dg}{\delta}\newc{\Dg}{\Delta}
\newc{\eg}{\varepsilon}
\newc{\zg}{\zeta}
\newc{\thg}{\theta}
\newc{\llg}{\lambda}\newc{\LLg}{\Lambda}
\newc{\kg}{\kappa}
\newc{\rg}{\rho}
\newc{\sg}{\sigma}\newc{\Sg}{\Sigma}
\newc{\tg}{\tau}
\newc{\fg}{\phi}\newc{\Fg}{\Phi}
\newc{\vfg}{\varphi}
\newc{\og}{\omega}\newc{\Og}{\Omega}
\newc{\pdh}{\partial}

\newc{\ccG}{{\cal G}}

\newc{\ii}[1]{\int_{#1}}
\newc{\iidx}[2]{{\dspl\int_{#1}~#2~dx}}
\newc{\bii}[1]{{\dspl \ii{#1} }}
\newc{\biii}[2]{{\dspl \iii{#1}{#2} }}
\newc{\su}[2]{\sum_{#1}^{#2}}
\newc{\bsu}[2]{{\dspl \su{#1}{#2} }}

\newc{\biiom}[1]{{\dspl\int_{\bdrom}~ #1 ~d\sg}}
\newc{\io}[1]{{\dspl\int_{\Og}~ #1 ~dx}}
\newc{\bio}[1]{{\dspl\int_{\bdrom}~ #1 ~d\sg}}
\newc{\bsir}{\bsu{i=1}{r}}
\newc{\bsim}{\bsu{i=1}{m}}

\newc{\iibr}[2]{\iidx{\bprw{#1}}{#2}}
\newc{\Intbr}[1]{\iibr{R}{#1}}
\newc{\intbr}[1]{\iibr{\rg}{#1}}
\newc{\intt}[3]{\int_{#1}^{#2}\int_\Og~#3~dxdt}
\newc{\itQ}[2]{\dspl{\int\hspace{-2.5mm}\int_{#1}~#2~dz}}
\newc{\mitQ}[2]{\dspl{\rule[1mm]{4mm}{.3mm}\hspace{-5.3mm}\int\hspace{-2.5mm}\int_{#1}~#2~dz}}
\newc{\mitQQ}[3]{\dspl{\rule[1mm]{4mm}{.3mm}\hspace{-5.3mm}\int\hspace{-2.5mm}\int_{#1}~#2~#3}}

\newc{\mitx}[2]{\dspl{\rule[1mm]{3mm}{.3mm}\hspace{-4mm}\int_{#1}~#2~dx}}
\newc{\mitmu}[2]{\dspl{\rule[1mm]{3mm}{.3mm}\hspace{-4mm}\int_{#1}~#2~d\mu}}
\newc{\iidmu}[2]{{\dspl\int_{#1}~#2~d\mu}}

\newc{\iidm}[3]{{\dspl\int_{#1}~#2~d #3}}

\newc{\itQmu}[2]{\dspl{\int\hspace{-2.5mm}\int_{#1}~#2~d\mu}}
\newc{\mitQmu}[2]{\dspl{\rule[1mm]{4mm}{.3mm}\hspace{-5.3mm}\int\hspace{-2.5mm}\int_{#1}~#2~d\mu}}

\newc{\mitQq}[2]{\dspl{\rule[1mm]{4mm}{.3mm}\hspace{-5.3mm}\int\hspace{-2.5mm}\int_{#1}~#2~d\bar{z}}}
\newc{\itQq}[2]{\dspl{\int\hspace{-2.5mm}\int_{#1}~#2~d\bar{z}}}

\newc{\pder}[2]{\dspl{\frac{\partial #1}{\partial #2}}}

\newc{\bdrom}{\bdry{\Og}}

\newc{\bilhom}{\mbox{Bil}(\mbox{Hom}(\RR^{nm},\RR^{nm}))}
\newc{\VV}[1]{{V(Q_{#1})}}

\newc{\ccA}{{\mathcal A}}
\newc{\ccB}{{\mathcal B}}
\newc{\ccC}{{\mathcal C}}
\newc{\ccD}{{\mathcal D}}
\newc{\ccE}{{\mathcal E}}
\newc{\ccH}{\mathcal{H}}
\newc{\ccF}{\mathcal{F}}
\newc{\ccI}{{\mathcal I}}
\newc{\ccJ}{{\mathcal J}}
\newc{\ccK}{{\mathcal K}}
\newc{\ccP}{{\mathcal P}}
\newc{\ccQ}{{\mathcal Q}}
\newc{\ccR}{{\mathcal R}}
\newc{\ccS}{{\mathcal S}}
\newc{\ccT}{{\mathcal T}}
\newc{\ccX}{{\mathcal X}}
\newc{\ccY}{{\mathcal Y}}
\newc{\ccZ}{{\mathcal Z}}

\newc{\bb}[1]{{\mathbf #1}}

\newc{\myprod}[1]{\langle #1 \rangle}
\newc{\mypar}[1]{\left( #1 \right)}
\newc{\BLLg}{\mathbf{\LLg}}

\newc{\mA}{\mathbf{A}}
\newc{\mB}{\mathbf{B}}
\newc{\mC}{\mathbf{C}}
\newc{\mD}{\mathbf{D}}
\newc{\mE}{\mathbf{E}}
\newc{\mF}{\mathbf{F}}
\newc{\mJ}{\mathbf{J}}
\newc{\mG}{\mathbf{G}}
\newc{\mP}{\mathbf{P}}
\newc{\mR}{\mathbf{R}}
\newc{\mQ}{\mathbf{Q}}
\newc{\mX}{\mathbf{X}}
\newc{\muu}{\mathbf{u}}
\newc{\mvv}{\mathbf{v}}

\newc{\mllg}{\mathbb{\lambda}}
\newc{\mLLg}{\mathbf{\LLg}}

\newc{\lspn}[2]{\mbox{$\| #1\|_{\Lsp{#2}}$}}
\newc{\Lpn}[2]{\mbox{$\| #1\|_{#2}$}}
\newc{\Hn}[1]{\mbox{$\| #1\|_{H^1(\Og)}$}}

\newc{\mynorm}[2]{\| #1\|_{#2}}

\newcommand{\RR}{{\rm I\kern -1.6pt{\rm R}}}

\newc{\itQQ}[2]{\dspl{\int_{#1}#2\,dz}}
\newc{\mmitQQ}[2]{\dspl{\rule[1mm]{4mm}{.3mm}\hspace{-4.3mm}\int_{#1}~#2~dz}}
\newc{\MmitQQ}[2]{\dspl{\rule[1mm]{4mm}{.3mm}\hspace{-4.3mm}\int_{#1}~#2~d\mu}}

\newc{\MUmitQQ}[3]{\dspl{\rule[1mm]{4mm}{.3mm}\hspace{-4.3mm}\int_{#1}~#2~d#3}}
\newc{\MUitQQ}[3]{\dspl{\int_{#1}~#2~d#3}}

\newc{\mccP}{\mathbb{P}}
\newc{\mccK}{\mathbb{K}}

\newc{\DKTmU}{\mccK(U)}
\newc{\DKTmUold}{(K_U(U)^{-1})^T}

\newc{\myPi}{\mathbf{W}}
\newc{\myIbar}{\bar{\ccI}_1}
\newc{\myIhat}{\hat{\ccI}_1}
\newc{\myIbreve}{\breve{\ccI}_0}

\newc{\mmk}{\mathbf{k}}

\newc{\mfu}{\mathbf{f_u}}
\newc{\mh}{\mathbf{h}}
\newcommand{\barrl}[2]{\barr{ll}\lefteqn{#1}\hspace{#2}&\\}

\begin{document}

\vspace*{-.8in}
\begin{center} {\LARGE\em Global Existence for some Cross Diffusion Systems with Equal Cross Diffusion/Reaction Rates.}

 \end{center}

\vspace{.1in}

\begin{center}

{\sc Dung Le}{\footnote {Department of Mathematics, University of
Texas at San
Antonio, One UTSA Circle, San Antonio, TX 78249. {\tt Email: Dung.Le@utsa.edu}
}}

\end{center}

\begin{abstract}
We consider some cross diffusion systems which is  inspired by models in mathematical biology/ecology, in particular the Shigesada-Kawasaki-Teramoto (SKT) model in population biology. We establish the global existence of strong solutions to systems for multiple species having equal either diffusion or reaction rates. The systems are given on bounded domains of arbitrary dimension.  \end{abstract}

\vspace{.2in}

\section{Introduction}\label{introsec}\eqnoset

In this paper, we study the global existence of following strongly coupled parabolic system of $m$ equations ($m\ge2$)  for the unknown vector $u=[u_i]_{i=1}^m$
\beqno{ep1}(u_i)_t=\Delta(u_ip_i(u))+u_ig_i(u),\quad (x,t)\in \Og\times(0,\infty).\eeq  
Here,  $p_i, g_i:\RR^m\to\RR$  are sufficienly smooth functions. Namely, $p_i\in C^2(\RR^m)$ and $g_i\in C(\RR^m)$. $\Og$ is a bounded domain with smooth boundary in $\RR^N$, $N\ge2$. 

The system is equipped with Dirichlet boundary and sufficiently smooth initial conditions
\beqno{ep1bc}\left\{\barr{l} \mbox{$u_i=0$  on $\partial \Og\times(0,\infty)$},\\ u_i(x,0)=u_{i,0}(x) ,\quad x\in \Og. \earr\right.\eeq

The consideration of \mref{ep1} is motivated by the extensively studied  model in population biology introduced by Shigesada {\it et al.} in \cite{SKT}
\beqno{e0}\left\{\barr{lll} u_t &=& \Delta(d_1u+\ag_{11}u^2+\ag_{12}uv)+k_1u+\bg_{11}u^2+\bg_{12}uv,\\v_t &=& \Delta(d_2v+\ag_{21}uv+\ag_{22}v^2)+k_2v+\bg_{21}uv+\bg_{22}v^2.\earr\right.\eeq   Here, $d_i,\ag_{ij},\bg_{ij}$ and $k_i$ are constants with $d_i>0$. Dirichlet or Neumann boundary conditions were usually assumed for \mref{e0}. This model was used to describe the population dynamics of {\em two} species densities  $u,v$ which move and interact with each other under the influence of their population pressures.

Of course, \mref{e0} is a special case of \mref{ep1} with $m=2$ and
$$p_i(u,v)=d_i+\ag_{i1}u+\ag_{i2}v,\; g_i(u,v)=k_i+\bg_{i1}u+\bg_{i2}v.$$
We will refer to the functions $p_i$'s (respectively, $g_i$'s) as the diffusion (respectively, raction) rates (see \cite{LM} for further discussions).

Under suitable assumptions on the constant parameters $\ag_{ij}$'s, $\bg_{ij}$'s and that $\Og$ is a planar domain ($N=2$), Yagi proved in \cite{yag} the global existence of (strong) positive solutions, with positive initial data. In this paper,  we will extend this investigation to multi-species versions of \mref{e0} for more than two species on bounded domains of arbitrary dimension $N$.

The global existence problem of \mref{ep1}, a fundamental problem in the theory of pdes. We can write \mref{ep1} in its general divergence form
\beqno{ep1div} u_t =\Div(A(u)Du)+f(u).\eeq
This a strongly coupled parabolic system with the diffusion matrix $A(u)$, the Jacobian of $[u_ip_i(u)]_1^m$, being a {\em full} matrix. We say that the system is weakly coupled if $A(u)$ is diagonal (i.e., $p_i$ depends only on $u_i$).

The key point in the proof of global existence of strong solutions of \mref{ep1div} is the a priori estimate of their spatial derivatives. In fact, it was established by Amann in \cite{Am2} that \mref{ep1} has a global strong solution $u$ if there is  some exponent $p>N$ such that for any $T\in (0,\infty)$
$$\limsup_{t\to T^-}\|Du\|_{L^p(\Og)}<\infty.$$

Thus, we need only prove that $\sup_{t\in(0,T) }\|Du\|_{L^p(\Og)}<\infty$ for all $T\in(0,\infty)$ and some $p>N$. With this a priori estimate, one can alternatively use the homotopy or fixed point approaches in \cite{letrans, dleANS,dlebook}, instead of semigroup theories in \cite{Am2}, to obtain the local/global existence of strong solutions.

The derivation of such estimates for \mref{ep1} is a difficult issue when $A(u)$ is full because the known techniques for scalar equations ($m=1$) are no longer applicable unless the matrix $A(u)$ are of special form, e.g., diagonal or triangular, these techniques can be partly applied together with some ad hoc arguments (see \cite{YW}). In this paper, we will consider \mref{ep1} with {\em full} diffusion matrix $A(u)$ of special forms where some nontrivial modifications of the classic methods can apply and yield new affirmative answers to the problem. 

Precisely, we study the case when either the diffusion or reaction rates are identical. Being inspired by the standard (SKT) system \mref{e0} where $p_i$ are a linear function in $u$, we consider a function $\Psi$ on $\RR$, a linear combination $L(u)$ of $u_i$'s, $L(u)=\sum_i a_iu_i$, and assume that for $i=1,\ldots,m$\beqno{pidef} p_i(u)=\llg_0+\Psi(L(u)).\eeq
We also assume that the reaction rates $g_i$'s satisfy the control growth $|g_i(u)|\le C+c_0\Psi(|L(u)|)$ for some positive constants $C,c_0$. We will establish the global existence of nonnegative strong solutions to \mref{ep1} with nonnegative initial data.

Clearly, \mref{e0} is the case when $d_1=d_2$, $\ag_{i1}=\ag_{j1}$, $\ag_{i2}=\ag_{j2}$ and $\Psi(s)=s$.

On the other hand, we can relax the assumption that the diffusion rates are identical as in \mref{pidef}. The trade off is that the reaction rates $g_i$'s are identical and satisfying the above control growth.

The paper is organized as follows. In \refsec{scalareqn}, we discuss some   regularity positivity  results  for strong solutions to scalar parabolic equations. Our main results on the system \mref{ep1} will be presented and proved in \refsec{equal}.

\section{Some facts on scalar equations} \label{scalareqn}\eqnoset In this section we consider the following scalar equation \beqno{scalarveqn}v_t = \Delta(P(v))+\Div(vb(v))+vg(v)\eeq in $Q=\Og\times(0,T)$ and
and study the smoothness, uniform boundedness and positivity of its {\em strong solution} $v$ under some special conditions on $P,g$ which will serve our purpose in discussing cross diffusion systems later.

To proceed, we first need the following parabolic Sobolev imbedding inequality.

\blemm{parasobolev} Let $r^*=p/N$ if $N>p$ and $r^*$ be any number in $(0,1)$ if $N\le p$. For any sufficiently nonegative smooth functions $g,G$ and any time interval $I$ there is a constant $C$ such that\beqno{paraSobo}\itQ{\Og\times I}{g^{r^*}G^p}\le C\sup_I\left(\iidx{\Og\times\{t\}}{g}\right)^{r^*}\itQ{\Og\times I}{(|DG|^p+G^p)}\eeq   If $G=0$ on the parabolic boundary $\partial\Og\times I$ then the integral of $G^p$ over $\Og\times I$ on the right hand side can be dropped.

Furthermore, if $r<r^*$ then for any $\eg>0$ we can find a constant $C(\eg)$ such that
\beqno{paraSobo1}\itQ{\Og\times I}{g^{r}G^p}\le C\sup_I\left(\iidx{\Og\times\{t\}}{g}\right)^{r}\itQ{\Og\times I}{(\eg|DG|^p+C(\eg)G^p)}\eeq
\elemm

\bproof For any $r\in(0,1)$ and $t\in I$ we have  via H\"older's inequality
\beqno{Sobo1}\iidx{\Og}{g^rG^p}\le \left(\iidx{\Og}{g}\right)^{r}\left(\iidx{\Og}{G^\frac{p}{1-r}}\right)^{1-r}.\eeq
If $r=r^*$ then $p/(1-r)=N_*=pN/(N-p)$, the Sobolev conjugate of $p$ if $N>p$ (the case $N\le p$ is obvious),  so that the Sobolev inequality gives
$$\left(\iidx{\Og}{G^\frac{p}{1-r}}\right)^{1-r}\le \iidx{\Og}{(|DG|^p+G^p)}.$$

Using the above in \mref{Sobo1} and integrating over $I$, we  easily obtain  \mref{paraSobo}. On the other hand, if $r<r^*$, then $p/(1-r)<N_*$. A simple contradiction argument and the compactness of the imbedding of $W^{1,p}(\Og)$ into $L^{p/(1-r)}(\Og)$ imply that for any $\eg>0$ there is $C(\eg)$ such that
$$\left(\iidx{\Og}{G^\frac{p}{1-r}}\right)^{1-r}\le \eg\iidx{\Og}{|DG|^p}+C(\eg)\iidx{\Og}{G^p}.$$ We then obtain \mref{paraSobo1}. \eproof

We now have the following a priori boundedness of solution of \mref{scalarveqn}.

\btheo{vboundthm} Consider a (weak or strong) solution $v$ to \mref{scalarveqn} in $Q=\Og\times(0,T)$. Assume that there are a function $\llg(v)$ and a number $\llg_0$ such that $\llg(v)\ge \llg_0>0$ and
\beqno{Qvcond}P_v(v)\ge \llg(v),\eeq \beqno{bvcond}|b(v)|\le g_1\llg(v),\eeq\beqno{fvcond}|g(v)|\le g_2\llg(v),\eeq
where $g_1,g_2$ are functions such that $g_1^2+g_2\in L^q(Q)$ for some $q>N/2+1$.

For $v\in \RR$ and $p\ge1$ consider the function \beqno{Fdef} F(v,p)=\int_0^v\llg^\frac12(s)s^{p-1}\,ds,\eeq and assume that \beqno{Fdefcond}|F(v,p)|\sim Cp\llg^\frac12(v)|v|^{p}\quad \mbox{for all $p$ and $v\in\RR$.}\eeq 

If $\|v\llg(v)\|_{L^1(Q)}$ is finite then $v,Dv$ are bounded and H\"older continuous in $\Og\times(\tau,T)$ for any $\tau\in(0,T)$. Their norms depend on $\|v\llg(v)\|_{L^1(Q)}$.

\etheo

The condition \mref{Fdefcond} is clearly verified if $\llg(v)$ has a polynomial growth in $|v|$.

\bproof We test the equation with $|v|^{2p-2}v$ and use integration by parts
$$\iidx{\Og}{\Delta(P(v))|v|^{2p-2}v}=-\iidx{\Og}{P_v(v)DvD(|v|^{2p-2}v)},$$
$$\iidx{\Og}{\Div(vb(v))|v|^{2p-2}v}=-\iidx{\Og}{vb(v)D(|v|^{2p-2}v)}.$$
Because $D(|v|^{2p-2}v)=(2p-1)|v|^{2p-2}Dv$ and the assumptions on $Q_v(v)$ and $b(v),g(v)$,
we easily get for {\em all} $p\ge1$ 
\beqno{vptest}\barrl{\sup_{(0,T)}\frac{1}{2p}\iidx{\Og}{|v|^{2p}}+(2p-1)\itQ{Q}{\llg(v)|v|^{2p-2}|Dv|^2}\le}{3cm} &C\itQ{Q}{g_1|\llg(v)||v|^{2p-1}|Dv|}+C\itQ{Q}{g_2|\llg(v)||v|^{2p}}.\earr\eeq
Applying Young's inequality $g_1|\llg(v)||v|^{2p-1}|Dv|\le \eg|v|^{2p-2}|Dv|^2+C(\eg)g_1^2|v|^{2p}$ for $\eg$ small,
$$\sup_{(0,T)}\frac{1}{2p}\iidx{\Og}{|v|^{2p}}+(2p-1)\itQ{Q}{\llg(v)|v|^{2p-2}|Dv|^2}\le C\itQ{Q}{(g_1^2+g_2)|\llg(v)||v|^{2p}}.$$

As $\llg(v)|v|^{2p-2}= F_v^2(v,p)$ by the definition \mref{Fdef}, for $g_3=g_1^2+g_2$ the above is 
$$\sup_{(0,T)}\frac{1}{2p}\iidx{\Og}{|v|^{2p}}+(2p-1)\itQ{Q}{|D(F(v,p)|^2}\le C\itQ{Q}{g_3|\llg(v)||v|^{2p}}.$$

Thus,  for $p\ge1$
$$\sup_{(0,T)}\iidx{\Og}{|v|^{2p}},\; \itQ{Q}{|D(F(v,p))|^2}\le Cp\itQ{Q}{g_3|\llg(v)||v|^{2p}}.$$

Applying the parabolic Sobolev inequality in \reflemm{parasobolev} with $g=|v|^p$ and $G=F(v,p)$, the above estimate  yields for $r=2/N$
$$ \left(\itQ{Q}{|v|^{2pr}|F(v,p)|^2}\right)^\frac{1}{1+r} \le Cp^{1+\frac{2}{1+r}}\itQ{Q}{g_3|\llg(v)||v|^{2p}}.$$
As $F(v)\sim Cp^{-1}\llg^\frac12(v)|v|^{2p}$ by \mref{Fdefcond}, we then obtain for $\cg=1+2/N$
$$ \left(\itQ{Q}{|v|^{2p\cg}\llg(v)}\right)^\frac{1}{\cg} \le Cp^{1+\frac{2}{\cg}}\itQ{Q}{g_3\llg(v)|v|^{2p}}.$$

H\"older's inequality yields
$$\itQ{Q}{g_3\llg(v)|v|^{2p}}\le C\left(\itQ{Q}{g_3^q\llg(v)}\right)^\frac{1}{q}\left(\itQ{Q}{|v|^{2pq'}\llg(v)}\right)^\frac{1}{q'}.$$
Let $d\mu=\llg(v)dz$. As we assume that  $g_1^2, g_2\in L^q(Q,d\mu)$,  $g_3\in L^q(Q,d\mu)$ and the first factor on the right hand side is finite. The above inequality is
\beqno{viterate}\|v\|_{L^{2p\cg}(Q,d\mu)} \le (2Cp)^{(1+\frac{2}{\cg})\frac{1}{2p}}\|v\|_{L^{2pq'}(Q,d\mu)}.\eeq

Because $q>N/2+1$, $q'<\cg=1+2/N$. Replacing $p$ by $pq'$ and defining $\cg_0=\cg/q'>1$). It follows that 
\beqno{viterate1}\|v\|_{L^{2p\cg_0}(Q,d\mu)} \le (2Cp)^{(\frac{1}{q'}+\frac{2}{\cg_0})\frac{1}{2p}}\|v\|_{L^{2p}(Q,d\mu)}.\eeq

Because $\cg_0>1$, we can apply the Moser iteration agument to show that $v$ is bounded. Indeed,  by taking $2p=\cg_0^i$ with $i=0,1,\ldots$. to the above estimate implies
$$\|v\|_{ L^{\cg^i}(Q,d\mu)}\le (2C)^{\cg_1}\cg^{\cg_2}\|v\|_{L^1(Q,d\mu)},$$ with $ \cg_1=(\frac{1}{q'}+\frac{1}{\cg_0})\sum_{i=0}^\infty \cg_0^{-i},\cg_2=(\frac{1}{q'}+\frac{1}{\cg_0})\sum_{i=0}^\infty i\cg_0^{-i}$. Letting $i\to\infty$ and using the fact that $\lim{p\to\infty\|v\|_{L^p(Q,d\mu)}=\|v\|_{L^\infty(Q,d\mu)}}$ (we will show that $d\mu$ is finite below) we obtain for some constant $C_0$
 that $\|v\|_{L^\infty(Q,d\mu)}\le C_0\|v\|_{L^1(Q,d\mu)}$.

As $\llg(v)$ is bounded below by a positive constant, this implies that   $v$ is bounded if $v\in L^1(Q,d\mu)$ is bounded. Furthermore, we now show that $d\mu$ is finite. Because $$ \itQ{|v|\ge1}{\llg(v)}\le \|v\|_{L^1(Q,d\mu)},$$ and $\llg(u)$ is bounded on the set $|v|<1$, we see that $d\mu$ is finite.

Once we show that $v$ is bounded, we obtain the {\em local} Harnack inequality (using both posive and negative power $p$ and cutoff functions) and so that $v$ is H\"older continuous. The argument is now classical and we refer the readers to the classical books \cite{LSU,Lieb} for details. It also follows that $Dv$ is bounded and H\"older continuous in $\Og\times(\tau,T)$ for any $\tau\in(0,T)$. Indeed, we can adapt the freezing coefficient method in \cite{GiaS} to establish this fact.  \eproof

\brem{vlarge} The conditions in the theorem and remarks need only hold only for $|v|$ large. This is easily to see if we make use of the cutoff function 
\beqno{vcut}\bar{v}_{(k)}=\left\{\barr{ll} v &\mbox{if $|v|\ge k$},\\k &\mbox{if $0<v< k$},\\
-k &\mbox{if $-k<v\le 0 $}\earr \right.\eeq
 with $k$ sufficiently large and observe that $D\bar{v}_k=0$ on the set $|v|<k$. \erem

\brem{vPsieqn} In connection with the systems considered in the next section, we consider the scalar equation
\beqno{veqn1} v_t= \llg_0 \Delta v +\Delta (\Psi(v)v)+ vg(v),\eeq
where $\llg_0>0$ and
$\Psi:\RR\to\RR$ be a $C^1$ function  and satisfying  for $|v|$ large
\beqno{Psicondv}\Psi(v),\;\Psi'(v)v\ge0.\eeq
Asume also that for $v\in \RR$ and $p\ge1$ the function \beqno{Fdef1} \hat{F}(v,p)=\int_0^v\Psi^\frac12(s)s^{p-1}\,ds\eeq satisfies \beqno{Fdefcond1}|\hat{F}(v,p)|\sim Cp\Psi^\frac12(v)|v|^{2p}\quad \mbox{for all $p$ and $v\in\RR$.}\eeq

This condition  allows us to apply \reftheo{vboundthm} with $P(v)=\llg_0 v +\Psi(v)v$ and $\llg(v)=\Psi(v)+\Psi'(v)v+\llg_0$. Thanks to \mref{Psicondv}, $\llg(v)$ satisfies \mref{Fdefcond}. Also, \mref{Fdef1} and \mref{Fdefcond1} imply that the function $F$ defined by \mref{Fdef} satisfies \mref{Fdefcond}.  We then apply \reftheo{vboundthm} to \mref{veqn1} and obtain that $v,Dv$ are bounded in $\Og\times(\tau,T)$ for any $\tau\in(0,T)$ and their norms are bounded in term of $\|v\|_{L^1(Q)}$ and $\|v\Psi(v)\|_{L^1(Q)})$.

We can also consider the scalar equation
\beqno{veqn2} v_t= \llg_0 \Delta v +\Delta (\Psi(|v|)v)+ vg(v),\eeq
and
$\Psi:\RR\to\RR$ be a $C^1$ function  and satisfying  for $v\ge0$ and large
\beqno{Psicondv1}\Psi(v),\;\Psi'(v)\ge0.\eeq
Indeed, we now define $\psi(v)=\Psi(|v|)$. We then have $\psi'(v)v=\Psi'(|v|)\mbox{sign}v v=\Psi'(|v|)|v|\ge0$ because of \mref{Psicondv1} $|v|\ge0$. Thus, $\psi$ satisfies \mref{Psicondv} and the theorem applies.

\erem

In applications we usually prefer that $v$ is nonnegative if the initial is. The following result serves this purpose.
\btheo{vposcoro} Let $a,g$ be $C^1$ functions on $\RR\times Q$ and $b$ be a bounded $C^1$ map from $Q$ into $\RR^N$. Assume that $a(w)\ge \llg_0$ for $w\ge0$ and $\llg_0$ is a positive constant. Also suppose that $a,g$ are bounded by a constant depending on $w$ in $(x,t)\in Q$.

Let $w$ be the strong solution to \beqno{veqn2a} \left\{\barr{ll}w_t= \Div(a(w,x,t)Dw)+\Div(wb) +wg(w,x,t), & \mbox{in $Q$}\\w(x,0)=w_0(x)&\mbox{on $\Og$}.
\earr \right.\eeq

If $w_0\ge0$ then $w\ge0$ on $Q$. \etheo

\bproof 
Because $w$ is a strong solution, there is a constant $M>0$ susch that $|w|\le M$.
We then truncate $a,g$ to $C^1$ function $\hat{a},\hat{g}$ which are constants for $v$ outside $[-M-1,M+1]$ and consider the equation \beqno{veqn2b} v_t= \Div(\hat{a}(|v|,x,t)Dv) + \Div(vb(x,t))+v\hat{g}(v,x,t),\eeq
with initial data $w_0$. 

We have $\hat{a}(|v|,x,t)\ge\llg_0$ and is bounded from above and $|v\hat{g}(v,x,t)|\le C|v|$ for some constant $C$. These facts and the classical theory of scalar parabolic equation show that \mref{veqn2} has a strong solution $v$. 

Let $v^+,v^-$ be the positive and negative parts of $v$. We test the equation with $v^-$. Using the facts that
$|v|=v^++v^-$, $|v|=v^-$ on the set $v^->0$, $v^+Dv^-=Dv^+Dv^-=0$ on the set $v^->0$,  we obtain
$$-\frac{d}{dt}\iidx{\Og}{(v^-)^2}-\iidx{\Og}{\hat{a}|Dv^-|^2}=\iidx{\Og}{[-bv^-Dv^-+(v^-)^2\hat{g}]}.$$

Because $b$ are bounded by a constant $C(M)$, applying Young's inequality $$\iidx{\Og}{|bv^-Dv^-|} \le \eg \iidx{\Og}{|Dv^-|^2} +C(\eg,M)\iidx{\Og}{(v^-)^2}.$$
Because $\hat{g}$ is bounded by a constant $C$  depending on $M$ and  $a(v)\ge\llg_0$, we can choose $\eg$ sufficiently small in the above inequality to arrive at
$$\frac{d}{dt}\iidx{\Og}{(v^-)^2}+\iidx{\Og}{|Dv^-|^2}\le C(M)\iidx{\Og}{(v^-)^2}.$$

Thus, we see that the function $$z(t)=\iidx{\Og}{(v^-)^2}$$ satisfies the differential inequality
$z'\le C_1z$ and $z(0)=0$ because the initial data $v_0\ge0$.  We then apply comparision theorem to the equation $y'=Cy$ with $y(0)=0$ which has the solution $y(t).=0$ We then have $z(t)=0$ for all $t\in(0,T)$. Hence, $v^-=0$ on $Q$ so that $v\ge0$. It follows that the solution $v$ of \mref{veqn2} also solves \mref{veqn2a}. By the uniqueness of strong solutions, $w=v\ge0$ in $Q$. \eproof

\section{Cross diffusion system with equal diffusion/reaction rates} \label{equal}\eqnoset

In this section, we consider the system \mref{ep1} and assume either that the diffusion rates $p_i$'s or the reaction rates are equal. We will always assume nonngative initial data $u_{i,0}$.

Throughout this section we will consider a nonnegative $C^1$ function $\Psi$ on $\RR$ satisfying \beqno{Psimaincond}\Psi'(s)\ge0 \mbox{ for $s\ge0$}.\eeq

\subsection{Equal diffusion rates:}
We first consider the following system of $m$ equations for $u=[u_i]_1^m$
\beqno{uWeqn}\left\{\barr{ll}(u_i)_t=\Delta(\llg_0 u_i +  \Psi(L(u))u_i)+u_ig_i(u)&\mbox{in $\Og\times(0,\infty)$},\\u_i(x,0)=u_{i,0}(x)&\mbox{on $\Og$},\earr\right.\eeq where $\llg_0>0$  and   $L(u)$ is a linear combination of $u_i$. That is, $ L(u)=\Sigma_{i=1}^m a_i u_i$ with $a_i>0$.
Assume that there are constants $ C_{ij},c_{ij}\ge0$ such that
\beqno{gicond} |g_i(u)|\le \sum_j (C_{ij}+c_{ij}\Psi(|u_i|)).\eeq

We have
\btheo{equaldiffthm} If $c_0=\max c_{ij}$ is sufficiently then \mref{uWeqn} has a unique nonnegative strong solution.\etheo

As we explained in the introduction, we need only establish a priopi the finiteness of $\sup_{(0,T_{max})}\|Du\|_{L^p(\Og)}$, with some $p>N$, for any strong solution $u=[u_i]_1^m$ of \mref{uWeqn} on $\Og\times(0,T_{max})$ for any $T_{max}\in(0,\infty)$. We will do this for $p=\infty$ via several lemmas.

\blemm{uposlemma} $u_i\ge0$ on $\Og\times(0,T_{max})$. \elemm
\bproof We can use \reftheo{vposcoro} to show first that $u_i\ge0$ on $Q=\Og\times[0,T]$ for any $0<T<T_{max}$ and all $i$. We rewrite the equation of $u_i$ as
\beqno{upsisys}\left\{\barr{ll}(u_i)_t=\Div(a_i(u_i,x,t)Du_i)+\Div(u_ib_i(x,t))+ u_ig_i(u) & \mbox{in $Q$},\\ u_i(x,0)=u_{i,0}(x)&\mbox{on $\Og$},\earr \right.\eeq where $$a_i(u_i,x,t)=\llg_0 + \Psi(L(u))+\partial_{u_i}\Psi(L(u))u_i,\; b_i(x,t)=\sum_{j\ne i} \partial_{u_j}\Psi(L(u(x,t))).$$

Following the proof of \reftheo{vposcoro}, because $u$ bounded on $Q$, $|L(u)|\le M$ for some constant $M$. We truncate the function $\Psi$ outside the interval $[-M-1,M+1]$ to obtain a bounded $C^1$ function $\psi$ satisfying: $\psi(s),\psi'(s)\ge0$ and $\psi(s)$ is a constant when $|s|\ge M+1$.

Denoting $\hat{v}=[|v_i|]_1^m$ for any vector $v=[v_i]_1^m$. We consider the system
\beqno{uvpsisys}\left\{\barr{ll}(v_i)_t=\Div(\hat{a}_i(v,x,t)Dv_i)+\Div(v_ib_i(x,t))+ v_ig_i(u) & \mbox{in $Q$}, \\v_i(x,0)=u_{i,0}(x)&\mbox{on $\Og$},\earr \right.\eeq  where $$\hat{a}_i(v,x,t)=\llg_0 + \psi(L(\hat{v}))+\partial_{v_i}\psi(L(\hat{v}))v.$$  Because $\psi'(s)\ge0$ for $s\ge0$ and $L(\hat{v})\ge0$, we have $\partial_{u_i}\psi(L(\hat{v}))v_i=\psi'(L(\hat{v}))a_i\mbox{sign}(v_i)v_i=\psi'(L(\hat{v}))a_i|v_i|\ge0$. We also have $\psi(L(\hat{v}))\ge 0$. Thus $\hat{a}_i(v,x,t)\ge\llg_0$ and bounded from above. The system \mref{uvpsisys} is a diagonal system with bounded continuous coefficients and has a unique strong solution $v$ according to the classical theory (e.g., see \cite[Chapter 7]{LSU}).

Applying the argument in the proof of \reftheo{vposcoro} to each equation in \mref{uvpsisys}, the system \mref{uvpsisys} has a nonnegative strong solution $v$, so that $\psi(\hat{v})=\psi(v)$, which also solves \mref{upsisys} by the definition of $\psi$, an extension of $\Psi$. By the uniqueness of strong solutions, $u_i=v_i\ge0$ in $Q$ for all $i$. \eproof

Next, define $W=L(u)$. The following lemma provides bounds of $W,DW$ that are independent of the number $M$, which was used only in establishing that $u_i\ge0$.
\blemm{Wbound} Let $W=L(u)\ge 0$.  Assume that
\beqno{Fpsidef} F(v,p):=\int_0^v\Psi^\frac12(s)s^{p-1}\,ds\sim Cp\Psi^\frac12(v)|v|^{2p}\quad \mbox{for all $p$ and $v\ge0$.}\eeq 

Then $W,DW$ are bounded in $\Og\times(\tau,T)$ for any $\tau\in(0,T)$ by a constant depending only on $\|W\|_{L^1(Q)},\|W\Psi(W)\|_{L^1(Q)}$.
\elemm

\bproof  
Taking a linear combination of the equations, we obtain
\beqno{Weqn1} W_t= \llg_0 \Delta W +\Delta (\Psi(W)W)+ f(u),\eeq 
where $f(u)=\sum_i a_i u_i g_i(u)$.
Because $u_i\ge0$ and $a_i>0$, $W$ is nonnegative and $|u_i|\le W$. Since $\Psi(s)$ is increasing for $s\ge0$, the assumption on $g_i$'s \mref{gicond} implies
$$|g_i(u)|\le \sum_j (C_{ij}+c_{ij}\Psi(|u_i|))\le \sum_j (C_{ij}+c_{ij}\Psi(W)).$$
Hence, $f$ satisfies for some positive constants $C$ and $c_0=\max c_{ij}$
\beqno{fgcond} |f(u)| \le C|W|(1+c_0\Psi(W)).\eeq

We then apply \reftheo{vboundthm} (to be precise, its \refrem{vPsieqn} and the equation \mref{veqn1}) with $v=W$, noting that $v=W\ge0$. The assumption \mref{Fpsidef} on $\Psi$ guarantees that \mref{Fdefcond1} is satisfied. We see that  the norms of $W,DW$ are bounded in $\Og\times(\tau,T)$ for any $\tau\in(0,T)$ by constants independent of $M$ but on $\|W\|_{L^1(Q)}$ and $\|W\Psi(W)\|_{L^1(Q)})$. The lemma follows.
\eproof

\brem{WL1bound}  If the constant  $c_0$ in \mref{fgcond} is sufficiently small then the norms  $\|W\|_{L^1(Q)}$ and $\|W\Psi(W)\|_{L^1(Q)})$ are bounded by a constant. Indeed, testing the equation of $W$ by $W$ and using \mref{fgcond}

\beqno{Wptest}\sup_{t\in (0,T)}\iidx{\Og\times\{t\}}{W^{2}}+\itQ{\Og\times(0,t)}{\Psi(W)|DW|^2}\le C\itQ{\Og\times(0,t)}{[1+c_0\Psi(W)]W^{2}}.\eeq

Applying the Sobolev inequality to the function $F(W,1)$ (see \mref{Fpsidef}) we find a constant $C(N)$ such that $$\iidx{\Og\times\{t\}}{\Psi(W)W^2}\le C(N)\iidx{\Og\times\{t\}}{\Psi(W)|DW|^2}.$$ Thus, using this, we see that if $c_0$ is sufficiently small then the integral of $Cc_0\Psi(W)W^2$ in  the inequality \mref{Wptest} can be absorbed to the left and we get
$$\sup_{t\in (0,T)}\iidx{\Og\times\{t\}}{W^{2}}+\itQ{\Og\times(0,t)}{\Psi(W)|DW|^2}\le C\itQ{\Og\times(0,t)}{W^{2}}.$$
This yields an integral Gr\"onwall inequality for $y(t)=\|W\|_{L^2(\Og\times\{t\})}$ on $(0,T)$ and shows that this norm is bounded by a universal constant on $(0,T)$. This fact and the  above inequality show that the left hand side quantities are bounded. We then make use of the parabolic Sobolev inequality to see that $\|W^{2\cg}\Psi(W)\|_{L^1(Q,d\mu)}$ is bounded by a constant. This implies $\|W\Psi(W)\|_{L^1(Q)})$ is bounded because $2\cg>1$.

\erem

{\bf Proof of \reftheo{equaldiffthm}:} 
We write the equation of $u_i$ in its divergence form
 $$(u_i)_t=\Div(aDu_i)+\Div(u_ib))+u_ig_i(u),$$ where $a=\llg +\Psi(W)$ and $b=D(\Psi(W))$.
 
Using the facts that $\Psi(W)\ge0$ (because $W\ge0$) and $W$ is bounded, we have $a\ge\llg_0$ and bounded from above. Also, $b=D(\Psi(W))$ are bounded. In addition,  $u_ig_i(u)$ is bounded because  $0\le u_i\le W/a_i$ which is bounded.   We then use 
the standard theory of scalar parabolic equation with bounded coefficients to show that $Du_i$ is bounded and H\"older continuous in $\Og\times(\tau,T)$ for any $\tau\in(0,T)$. \eproof

\subsection{Equal reaction rates:}

We now present two examples which relax the assumption of equal diffusion rates $p_i$'s. However, we have to consider equal reaction rates $g_i$'s and restrict ourselves to the case of systems of two equations.

In the sequel, we will always assume that $\Psi$ is a $C^1$ function on $\RR$ such that
\beqno{PsiWpos}\Psi(s), \Psi'(s)\ge 0 \mbox{ and }  \Psi(s)\ge s \mbox{ for $s\ge0$}.\eeq

We consider first the following system
\beqno{YW}\left\{\barr{lll} u_t&=&\Delta(\llg_0 u +u\Psi(L(u,v)) +\eg_0  a\Delta(uv)+ug(u,v),\\
v_t&=&\Delta(\llg_0 v + v\Psi(L(u,v))-\eg_0  b\Delta(uv)+vg(u,v).\earr\right.\eeq
Here, $L(u,v)=bu+av$. $\llg_0,\eg_0,a,b$ are positive constants.
Regarding the reaction term, we also assume that there are positive constants $C,c_0$ such that (compare with \mref{gicond})
\beqno{guvcond} |g(u,v)|\le C+c_0\Psi(|L(u,v)|) \mbox{ for all $u,v\in\RR$}.\eeq
 
We consider nonnegative initial data $u_0, v_0$ for $u,v$. 

\btheo{uvthm}If $\eg_0, c_0$ are sufficiently small then the system \mref{YW} has a unique global strong solution $(u,v)$ with $u,v\ge0$. \etheo

We need the following proposition which will be useful later.
\bprop{uvprop}
We consider a strong solution $(u,v)$ with nonnegative initial data $u_0,v_0$ to the following system
\beqno{YW0}\left\{\barr{lll} u_t&=&\Delta(\llg_0 u +u\Psi(L(u,v)) +\eg_0  a\Delta(u|v|)+ug(u,v),\\
v_t&=&\Delta(\llg_0 v + v\Psi(L(u,v))-\eg_0  b\Delta(u|v|)+vg(u,v).\earr\right.\eeq

For any $\eg_0>0$ we have that $u,v$ and $Du,Dv$ are bounded. Also $u\ge0$ in $Q$.

If $\eg_0$ are sufficiently small then $v$ is also nonnegative in $Q$.
\eprop

\bproof The proof will be divided into several steps. First of all, taking a linear combination of the above two equations,  we see that $W=L(u,v)$ satisfying
\beqno{Weqn1a} W_t= \llg_0 \Delta W +\Delta (\Psi(W)W)+ Wg(u,v).\eeq

{\bf Step 1:} We show that $W,DW$ are bounded and $W\ge0$.

For a given strong solution $(u,v)$ of \mref{YW0} we consider the the equation
\beqno{Weqn1b} w_t= \llg_0 \Delta w +\Delta (\Psi(|w|)w)+ wg(u,v)\eeq and the initial data $w_0=au_0+bv_0\ge0$.
We proved in \reftheo{vboundthm} that this equation has a strong solution $w$  and, by \reftheo{vposcoro}, $w\ge0$. By uniqueness of strong solutions, $W=w$ so that $W\ge0$.

Now, from the proof of \reflemm{Wbound} we see that $W,DW$ are bounded in $\Og\times(\tau,T)$ for any $\tau\in(0,T)$ in terms of $\|W\|_{L^1(Q)}$ and $\|W\Psi(W)\|_{L^1(Q)}$ The latter two norms can be bounded by a constant if $c_0$ is sufficiently small (see \refrem{WL1bound}). 

We should note that because we already proved that $W\ge0$, hence we do not need here the fact that $u,v\ge0$ (which will be established later) as before in \reflemm{Wbound} but the conditions $\Psi(s),\Psi'(s)\ge0$ for $s\ge0$ in \mref{PsiWpos} and that $|g(u,v)|\le C+c_0\Psi(|W|)$ in \mref{guvcond} (see equation \mref{veqn2} of \refrem{vPsieqn}).

{\bf Step 2:} We prove that $u\ge0$. We write the equation of $u$ in its divergence form
\beqno{udiveqn}u_t=\Div(ADu)+\Div(uB)+ug(u,v)\eeq
with $A=\llg_0 +\Psi(W)+\eg_0a|v|$, $B=-\Psi'(W)DW+\eg_0 aD(|v|)$.

Again, we can assume that $u,v$ are locally bounded as in the proof of \reflemm{uposlemma}. Because $W,DW$ are bounded, we apply \reftheo{vposcoro} to prove that $u\ge0$. 

{\bf Step 3:} We now prove that $u$ is bounded by using the iteration argument in \reftheo{vboundthm}. 

We multiply the above equation \mref{udiveqn} by $u^{2p-1}$, recall that $u\ge0$, and follows the proof of \reftheo{vboundthm} to get
\beqno{uiter}\frac{d}{dt}\iidx{\Og}{u^{2p}}+(2p-1)\iidx{\Og}{Au^{2p-2}|Du|^2}\le \iidx{\Og}{\Div(uB)u^{2p-1}}+\iidx{\Og}{g(u,v)u^{2p}}\eeq

From the definition of $B$ we need to study the following two terms on the right of \mref{uiter} 
\beqno{adiv}-\iidx{\Og}{\Div(u\Psi'(W)D(W))u^{2p-1}},\;\iidx{\Og}{a\Div(uD(|v|))u^{2p-1}}.\eeq

The first one can be treated easily, using the fact that $W$ is bounded (see also below). We consider the second term. We have
$$\iidx{\Og}{a\Div(uD(|v|))u^{2p-1}}
=-(2p-1)\iidx{\Og}{auD(|v|)u^{2p-2}Du}.$$

For each $t>0$ we split $\Og=\Og^+(t)\cup\Og^-(t)$ where $\Og^+(t)=\{x:\,v(x,t)\ge0\}$. 
Since $aDv=DW-bDu$ and on $\Og^+(t)$, $D(|v|)=Dv$, we have that the integral over $\Og^+(t)$ of $-auD(|v|)u^{2p-2}Du$ is
$$ \iidx{\Og^+(t)}{u^{2p-1}(-DW+bDu)Du}=-\iidx{\Og^+(t)}{u^{2p-1}DWDu}+\iidx{\Og^+(t)}{bu^{2p-1}|Du|^2}.$$

Because $DW$ is bounded in $\Og\times(\tau,T)$ for any $\tau\in(0,t)$, it follows that for any $\eg>0$ there is $c_1(\eg)$ such that
$$\iidx{\Og}{u^{2p-1}|DWDu|}\le \iidx{\Og}{(\eg u^{2p-2}|Du|^2 +c_1(\eg)u^{2p})}.$$
Choosing $\eg$ small, the integral of $u^{2p-2}|Du|^2$ can be absorbed to the integral of $\llg_0 u^{2p-2}|Du|^2$ in the left of \mref{uiter}. This argument also applies to the first integral in \mref{adiv}.

Meanwhile, on the set $v\ge0$, as $W\ge bu\ge0$  so that $\Psi(W)\ge W\ge bu$ (by the assumption \mref{PsiWpos} on $\Psi$). Thus, the integral over $\Og^+(t)$ of $bu^{2p-1}|Du|^2$ can also be absorbed to the integral over $\Og^+(t)$ of $\Psi(W)u^{2p-2}|Du|^2$ in $Au^{2p-2}|Du|^2$ of the left of \mref{uiter}.

On $\Og^-(t)$, $D(|v|)=-Dv$, we have that the integral over $\Og^-(t)$ of $-auD(|v|)u^{2p-2}Du$ is
$$ \iidx{\Og^-(t)}{u^{2p-1}(DW-bDu)Du}=\iidx{\Og^-(t)}{u^{2p-1}DWDu}-\iidx{\Og^-(t)}{bu^{2p-1}|Du|^2}.$$
The first integral on the right hand side can be handled as before. The second integral is nonnegative and can be dropped.

Putting these togetter, we then obtain for all $p\ge1$
$$\frac{d}{dt}\iidx{\Og}{u^{2p}}+\iidx{\Og}{u^{2p-2}|Du|^2}\le C\iidx{\Og}{u^{2p}}.$$ This allows to obtain a bound for $\|u\|_{L^\infty(Q)}$ in terms of $\|u\|_{L^1(Q)}$ (see \reftheo{vboundthm}). Let $p=1$ in the above inequality to get a Gr\"onwall inequality for $\|u\|_{L^2(\Og)}^2$. We see that $\|u\|_{L^2(\Og)}^2$, so is $\|u\|_{L^1(\Og)}$, is bounded on $(0,T)$.

Once we prove that $u$ and $W,DW$ are  bounded we then use a cutoff function and repeat a similar argument to the above one in order to obtain local strong/weak Harnack inequalities. It follows that $u$ is H\"older continuos. This is a standard procedure and the readers are referred to the book \cite{Lieb}. It also follows that $Du$ is bounded.

{\bf Step 4:} We show that $v$ is bounded. This is easy because $v=(W-bu)/a$ and $u,Du$ and $W,DW$ are bounded. We should note that in the above steps we have not imposed any assumptions on $\eg_0$. Thus, the first assertion of the Proposition was proved.

{\bf Step 5:} Finally, we prove that $v\ge0$. First of all, we write the equation of $v$ in its divergence form
$$ v_t=\Div(A_1Dv)+\Div(|v|B_1)+vg,$$ where $A_1=\llg_0 +\Psi(W)-\eg_0bu\mbox{sign}(v)$, $B_1=\Psi'(W)DW+\eg_0 bD(u)$.

Since $u$ is bounded by a constant independent of $\eg_0$ and $\Psi(W)\ge0$, we can choose $\eg_0$ small such that $A_1\ge\llg_0/2$. Also, $B_1$ is bounded because $W,DW$ and $Du$ are. The proof of \reftheo{vposcoro} applies and shows that $v\ge0$. \eproof

{\bf Proof of \reftheo{uvthm}:} From \refprop{uvprop} the system \mref{YW0} has a strong solution $(u,v)$ which also solves \mref{YW}. By uniqueness of strong solutions, we see that strong solution $(u,v)$, and its spatial derivatives, of \mref{YW} are bounded uniformly in terms of the data. Because $\|Du\|_{L^\infty(\Og)},\|Dv\|_{L^\infty(\Og)}$ do not  blow up in any time interval $(0,T)$, the solution exists globally. \eproof

We also consider the following system
\beqno{YWz}\left\{\barr{lll} u_t&=&\Delta(\llg_0 u +u\Psi(L(u,v)) +\eg_0  a\Delta(uv)+ug(u,v),\\
v_t&=&\Delta(\llg_0 v + v\Psi(L(u,v))+\eg_0  b\Delta(uv)+vg(u,v).\earr\right.\eeq
Here, $L(u,v)=bu-av$ and $\eg_0,a,b$ are positive constants.

We then have the following result similar to \reftheo{uvthm} without the assumption on the smallness of $\eg_0$. However, we have to strengthen the condition \mref{PsiWpos} by assuming in addition that
\beqno{PsiWpos1} \Psi(s)\ge0 \quad \forall s\in\RR.\eeq

\btheo{uvthmz}If $c_0$ in the assumption \mref{guvcond} is sufficiently small then the system \mref{YW} has a unique global strong solution $(u,v)$ with $u,v\ge0$. \etheo

\bproof Following the proof of \reftheo{uvthm}, we consider a strong solution $(u,v)$ with the same initial data to the following system
\beqno{YWz0}\left\{\barr{lll} u_t&=&\Delta(\llg_0 u +u\Psi(L(u,v)) +\eg_0  a\Delta(u|v|)+ug(u,v),\\
v_t&=&\Delta(\llg_0 v + v\Psi(L(u,v))+\eg_0  b\Delta(u|v|)+vg(u,v).\earr\right.\eeq

For any $\eg_0>0$ we will prove that $u,v$ and $Du,Dv$ are bounded. We also show that $u,v\ge0$ in $Q$. We follow the proof of \refprop{uvprop} and provide necessary modifications.

Let $W=bu-av$. Taking a linear combination of the two equations, we can follows Step 1 of the proof of \refprop{uvprop} to show that $W,DW$ are bounded. Note that we cannot prove that $W\ge0$ as before because its initial data $bu_0-av_0$ is not nonnegative. 

Similarly, Step 2  also yields that $u\ge0$. We need to change the argument in Step 3 of the proof to prove that $u,Du$ are bounded. We test the equation of $u$ by $u^{2p-1}$. As in Step 3, we need to consider the following term on the right hand side of \mref{uiter}
$$\iidx{\Og}{a\Div(uD(|v|))u^{2p-1}}=-(2p-1)\iidx{\Og}{auD(|v|)u^{2p-2}Du}.$$

We again split $\Og=\Og^+\cup\Og^-$ where $\Og^+=\{v\ge0\}$. Because $av=bu-W$ (instead of $av=W-bu$ as before) we need to interchange $\Og^+,\Og^-$ in the previous argument. Namely, the integral over $\Og^+$ now contributes a nonnegative term to the left and an integral of $u^{2p}$ to the right. Meanwhile, on $\Og^-$ we have $W =bu-av\ge bu\ge0$ so that $\Psi(W) \ge bu$ and the integral over $\Og^-$ of $bu^{2p-1}|Du|^2$ now can be absorbed to the left hand side. The proof then continues to prove that $u,Du$ are bounded.

Using $v=(bu-W)/a$, we see that $v,Dv$ are bounded.

We now show that $v\ge0$, without the assumption that $\eg_0$ is small. We slightly modify Step 5 of \refcoro{uvprop}. We write the equation of $v$ as
$$ v_t=\Div(A_2Dv)+\Div(\eg_0uD(|v|))+\Div(vB_2)
+ vg(u,v).$$
Here, $A_2=\llg_0 +\Psi(W)$, $B_2=\eg_0\mbox{sign}(v)Du+D\Psi(W)$. We follow the proof of \reftheo{vposcoro} and test the equation with $v^-$. We need to consider the integral of $\Div(\eg_0uD(|v|))v^-$ on the right hand side. Using integration by parts and the fact that $D(|v|)=Dv^++Dv^-$,
$$\iidx{\Og}{\Div(\eg_0uD(|v|))v^-}=-\iidx{\Og}{\eg_0uD(|v|)Dv^-}=-\iidx{\Og}{\eg_0u|Dv^-|^2}.$$
Because $u\ge0$, the last term provides a nonnegative term on the left hand side. Meanwhile, we have that $A_2\ge \llg_0$ and $A_2,B_2$ are bounded (as $u,Du,W,DW$ are bounded). 
We obtain as in the proof of \reftheo{vposcoro} a Gr\"onwall inequality of $\|v^-\|_{L^2(\Og)}$ and conclude that $v^-=0$ on $Q$. Thus, $v$ is nonnegative. The proof is complete.
\eproof

\bibliographystyle{plain}

\end{document}